\documentclass[a4paper]{article}

\usepackage[latin1]{inputenc}
\usepackage{amsmath}
\usepackage{rotating}
\usepackage{calc}
\usepackage{graphics}
\usepackage{color}

\usepackage[pdftitle={Monodromy calculatons of fourth order equations
  of Calabi-Yau type}, pdfauthor={Christian van Enckevort and Duco
  van Straten},ps2pdf]{hyperref}
\newcommand{\hlink}[1]{\href{#1}{\texttt{}#1}}
\newcommand{\dburl}{http://enriques.mathematik.uni-mainz.de/enckevort/db}
\newcommand{\eqnno}[2][4]{\href{\dburl/summary.php?order=#1\&eqnno=#2}{#2}}

\usepackage{array}
\newcolumntype{C}{>{$}c<{$}}
\newcolumntype{L}{>{$}l<{$}}
\newcolumntype{R}{>{$}r<{$}}

\usepackage{myref}
\usepackage{mysymb}

\newcommand{\Id}{\mathrm{Id}}
\newcommand{\Grass}{\mathrm{Grass}}
\newcommand{\LGrass}{\mathrm{LGrass}}
\renewcommand{\P}{\mathbb{P}}
\renewcommand{\L}{{\mathbb{L}}}
\newcommand{\lra}{{\longrightarrow}}
\newcommand{\Auteq}{\mathop{\mathrm{Auteq}}}

\newtheorem{hyp}{Hypothesis}

\begin{document}
\title{Monodromy calculations of fourth order equations of Calabi-Yau
  type}
\author{Christian van Enckevort and Duco van Straten}
\maketitle

{\renewcommand{\thefootnote}{}}

\begin{abstract}
  This paper contains a preliminary study of the monodromy of certain
  fourth order differential equations, that were called of Calabi-Yau
  type in \cite{AZ}.  Some of these equations can be interpreted as
  the Picard-Fuchs equations of a Calabi-Yau manifold with one complex
  modulus, which links up the observed integrality to the conjectured
  integrality of the Gopakumar-Vafa invariants. A natural question is
  if in the other cases such a geometrical interpretation is also
  possible. Our investigations of the monodromies are intended as a
  first step in answering this question.  We use a numerical approach
  combined with some ideas from homological mirror symmetry to
  determine the monodromy for some further one-parameter models.
  Furthermore, we present a conjectural identification of the
  Picard-Fuchs equation for 5 new examples from Borceas list and
  conjecture the existence of some new Calabi-Yau three folds. The
  paper does not contain any theorems or proofs but is, we think,
  nevertheless of interest.
\end{abstract}

\section{Introduction} 
A differential operator of order $n$ on $\P^1$ has the form 
\begin{equation}
\label{eq:L} 
  L:= a_n(z) \frac{d^n}{dz^n} + a_{n-1}(z)\frac{d^{n-1}}{dz^{n-1}} +
  \dots + a_0(z),
\end{equation}
where the $a_i(z)$ are polynomials. The set $\Sigma \subset \P^1$ of
singular points is given by the zeros of $a_n(z)$ and possibly
$z=\infty$.  The solutions to the equation $L y=0$ can be considered
as a $\C$-local system $\L$ of rank $n$ on $S:=\P^1\setminus \Sigma$.
After the choice of a base point $s \in \P^1 \setminus \Sigma$, the
information of $\L$ is given by the monodromy representation
\[
  \pi_1(S,s) \lra \Aut(\L_s) = \mathrm{Gl}_n(\C)
\]
A power series $y_0(x) \in \Z[[x]]$ that satisfies a homogeneous
linear differential equation as above is a G-function and a folklore
conjecture that goes back to Bombieri and Dwork states that all such
power series and differential operators have a \emph{geometrical
  origin} (see \cite{KZ}).  This means that the operator should occur
as a \emph{factor} of a Picard-Fuchs operator describing the variation
of a cohomology of a family $\rho: \mathcal{Y} \lra \P^1$, with
singular fibres over $\Sigma$ and defined over a number field. The
local system $\L$ should then be a summand of a local system
$\L_{\C}:=R^d\rho_*(\C_{\mathcal{Y}})_{|S}$, where $d$ is the complex
dimension of the fibres of $\rho$.  It follows among other things that
the equation has regular singularities with all exponents rational. It
can be shown that the set of power series of geometric origin in this
sense is closed under the ordinary product of power series and under
the coefficientwise Hadamard product of series.  On the level of local
systems, the Cauchy product corresponds to the tensor product, whereas
the Hadamard product correspond to the convolution of local systems.
We refer to the books \cite{An} and \cite{Ka} for details.

The fourth order equations in \cite{AZ} were collected with a stricter
notion of geometrical origin in mind: by requiring that the operator
admits an invariant symplectic form and gives rise to integral
instanton numbers, it starts making sense asking for the existence of
a one-parameter family $\mathcal{Y} \lra \P^1$ of \emph{Calabi-Yau
  three folds}, whose associated Picard-Fuchs operator for
$\Hgrp^3(Y_s)$ is the given one. The instanton numbers then should
have the interpretation of counting curves on a mirror manifold $X$
with Picard number one.  The first $14$ equations in the list are in
fact the much studied hypergeometric cases (see \cite{CdOGP},
\cite{M1}, \cite{KT}, \cite{BvS}, \cite{V}, \cite{DM}). Mirror pairs of
Calabi-Yau threefolds obtained from Batyrev's polar duality of
reflexive polytopes \cite{Ba} yield a plethora of examples but usually
with high Picard number (see \cite{CYH}).  By taking restrictions to
carefully chosen one-dimensional sub-loci these examples sometimes
give rise to equations of Calabi-Yau type, but the instanton numbers
computed in this way represent \emph{sums over different homology
  classes} and there will not exist a Calabi-Yau three fold $X$ with
Picard number one with the given instanton numbers. Case $15$ is an
example of this phenomenon: it is the equation belonging to the
diagonal restriction of Calabi-Yau family in $\P^3 \times \P^3$ (see
\cite{BvS}). The list contains many more of such examples. The question
is how can one see this from the differential equation alone.

In order to find the cases that are potentially of strict geometric
origin, we remark that a geometrical local system $\L_{\C}$ carries a
integral lattice $\L_{\Z} =R^d\rho_*(\Z_{\mathcal{Y}})|_S$ and that
Poincaré-duality provides it with a \emph{unimodular} pairing
$\langle\cdot,\cdot\rangle$, which in our case is alternating. Hence
the monodromy representation is in the symplectic group $\Sp(4,\Z)$.
For differential equations of hypergeometric type, the monodromy
representation is explicitly known, essentially because the associated
local system is rigid (Levelt's theorem, \cite{BH},\cite{Ka}).  This
leads to the $14$ hypergeometric cases mentioned above.  For equations
with three singular points which are not of hypergeometric type or for
equations with more than three singular points, the monodromy
representation is in general not determined by local data alone and we
have the problem of accesory parameters. We do not know of any general
method to determine the monodromy representation in such cases.  We
use a brute force numerical approach combined with ideas from
homological mirror symmetry to conjecturally determine the monodromy
for some further one-parameter models.

\textbf{Acknowledgements}: We would like to thank G.~Almkvist,
C.~Doran, A.~Klemm, and W.~Zudilin for their interest in the project.
In particular we thank C.~Doran for his explanation of the integral
basis and A.~Klemm for the suggestion of using the genus one instanton
numbers as an extra integrality check and for explanations on higher
genus computations.

\section{Sketch of Homological Mirror Symmetry}
According to Kontsevich \cite{Ko1}, the phenomenon of mirror symmetry
between Calabi-Yau spaces $X$ and $Y$ should be formulated in terms of
equivalence of categories. To a Calabi-Yau space $X$ one can associate
two triangulated categories, namely the derived category of coherent
sheaves $D^b(X)$ and a derived Fukaya-category $D\mathcal{F}(X)$ of
lagrangian cycles (graded, with local systems on them) in $X$
(see~\cite{FOOO}).  The first category depends only on the holomorphic
moduli, the second only on the symplectic (or Kähler) moduli.
Mirror symmetry between Calabi-Yau spaces $X$ and $Y$ is then
expressed as equivalences of categories.
\[
  \mathrm{Mir}: D^b(X) \stackrel{\approx}{\lra}
  D\mathcal{F}(Y),\quad D^b(Y) \stackrel{\approx}{\lra}
  D\mathcal{F}(X)
\] 
These equivalences induce isomorphisms between the corresponding
$K$-groups. Via the Chern character they descend to cohomology:
\[
  \mathrm{mir}: \Hgrp^{\text{ev}}(X,\Q) \stackrel{\approx}{\lra}
  \Hgrp^d(Y,\Q),\quad \Hgrp^{\text{ev}}(Y,\Q) 
  \stackrel{\approx}{\lra} \Hgrp^d(X,\Q), 
\]
where $d = \dim_C X = \dim_C Y$. This also induces an isomorphism
between the Kähler moduli $\Hgrp^{1,1}(X)$ of $X$ and the complex
moduli of $\Hgrp^{d-1,1}(Y)$ of $Y$.  In the Strominger-Yau-Zaslow
picture of mirror symmetry (see \cite{SYZ}, \cite{Gr}) $X$ and $Y$ are
represented as (real) singular torus fibration over a common base $B$.
The fibres are dual tori and mirror symmetry should correspond to
fibrewise T-duality.  From this one can get some intuitive
understanding of the mirror transformation on objects. In particular,
the structure sheaf $\mathcal{O}_p$ of a point $p \in X$ gets mapped
to a SYZ-fibre $\mathbf{T}$ (with a local system on it) in $Y$ and the
structure sheaf $\mathcal{O}_X$ should map to the image $\mathbf{S}$
of a section $\sigma: B \lra Y$ of the fibration.

For any pair $(\mathcal{E},\mathcal{F})$ of objects of $D^b(X)$ the Euler 
bilinear form is defined by
\[
  \langle\mathcal{E}, \mathcal{F}\rangle := \chi(\mathcal{E},\mathcal{F}) = 
  \sum_i (-1)^i \dim \Hom(\mathcal{E}, \mathcal{F}[i]).
\]
which by Serre duality and triviality of the canonical bundle is
$(-1)^d$ symmetric.  It descends via the Chern-character to a bilinear
form $\langle \cdot, \cdot \rangle$ on the cohomology
$\Hgrp^{\text{ev}}(X,\Q)$ of $X$, which by Riemann-Roch is given by
\[
  \langle \alpha,\beta \rangle = \int_X \tilde{\alpha} \cup \beta \cup \td(X),
\]
where $\tilde{\alpha} = (-1)^k \alpha$ for $\alpha \in
\Hgrp^{2k}(X,\Q)$. 

Under the mirror transformation the form $\langle \cdot, \cdot \rangle$ should
correspond to the intersection form $\langle \cdot, \cdot \rangle$ of the
corresponding lagrangians. One instance of this can easily be checked
\[ 
  \langle\mathcal{O}_p, \mathcal{O}_X\rangle = 1 
  = \langle \mathbf{T}, \mathbf{S} \rangle.
\]

\section{Monodromy in  one-parameter models}
From now on we assume that $X$ and $Y$ are strict Calabi-Yau
three-folds and furthermore that they satisfy $h^{2,1}(Y) = 1 =
h^{1,1}(X)$. This is the case of so called one-parameter models: $Y$
varies in a one-dimensional moduli space and $X$ has one Kähler
modulus, i.e., $\Pic(X)=\Z$.  In such a case one has $\dim \Hgrp^3(Y)
= 4 = \dim \Hgrp^{\text{ev}}(X)$.

To be specific, we assume that we have a proper map $\rho: \mathcal{Y}
\lra \P^1$, smooth outside singular fibres that sit over points from
$\Sigma \subset \P^1$ and furthermore that $Y$ is the fibre over a
base-point $s \in \P^1 \setminus \Sigma=:S$.  As the geometrical
monodromy along a path $\gamma \in \pi_1(S,s)$ can be realised as a
symplectic map $M(\gamma): Y_s \lra Y_s$, $M(\gamma)$ induces an
autoequivalence of its symplectic invariant $D\mathcal{F}(Y)$, thus
setting up a homomorphism
\[
  \pi_1(S,s) \lra \Auteq(D\mathcal{F}(Y))
\]
which is a refined version of the ordinary monodromy representation of
$\pi_1(S,s)$ on $\Hgrp^{\text{odd}}(Y)$.  The group $\pi_1(S,s)$ is
generated by paths that encircle one of the singular fibres of the
family. The induced transformation is determined by the specific
properties of the singular fibre.  If the fibre aquires the simplest
type of singularity, namely an $A_1$-singularity (`conifold'), there
is a vanishing lagrangian $3$-sphere. The geometrical monodromy is
then a Dehn-twist along this sphere and its effect on homology is given
by the classical Picard-Lefschetz transformation \cite{Le},
\cite{AGV}, \cite{Lo}:
\[
  \alpha \mapsto S_{\delta}(\alpha) := \alpha - \langle \delta, \alpha
  \rangle \delta 
\] 
where $\delta$ is the homology class of the vanishing cycle. 
In the situation of mirror symmetry there also will be a point of
degeneration with maximal unipotent monodromy. The fibre will typically
have normal crossing singularities and there will be a `vanishing
$n$-torus', invariant under the monodromy.

Using the mirror equivalence $\mathrm{Mir}$ we get a representation
\[
  \pi_1(S,s) \lra \Auteq(D^b(X))
\]
and one may ask what sort of autoequivalences correspond to
specific types of degenerations of $Y$.

In \cite{ST00} Seidel and Thomas described a type of autoequivalence
in $D^b(X)$ to mirror a symplectic Dehn-twist.  It is the
Seidel-Thomas twist $T_{\mathcal{E}}$ by a so called spherical object
$\mathcal{E}$ of $D^b(X)$, which has the property that
$\dim(\textup{Ext}^*(\mathcal{E},\mathcal{E})) = \dim
\Hgrp^*(\mathbf{S})$ and is given by the triangle
\[ 
\lra (\mathcal{E}, \mathcal{F})\otimes \mathcal{E} \lra \mathcal{F} \lra
  T_{\mathcal{E}}(\mathcal{F})\stackrel{+1}{\lra}
\]
The structure sheaf $\mathcal{O}_X$ is the basic spherical object in
$D^b(X)$, but also each line bundle $L \in \Pic(X)$ is spherical.
Another particularly simple type of autoequivalence is the operation
$\otimes L$ of tensoring with a line bundle $L$. Note that
$\mathcal{O}_p \otimes L = \mathcal{O}_p$. This fits on the mirror
side to the monodromy tranformation around a point of maximal
unipotent monodromy, with invariant vanishing torus $\mathbf{T}$.

Let us write out these transformations on the level of cohomology.
Let $L = \mathcal{O}(H)$ be the ample generator of $\Pic(X)$.
The powers  $1,H,H^2,H^3$ form a basis for $\Hgrp^{\text{ev}}(X,\Q)$.
With respect to this basis, the matrix $T$ of tensoring with
$L$ is given by
\begin{equation}
\label{eq:DM0}
  T = \begin{pmatrix}
          1 & 0 & 0 & 0 \\
          1 & 1 & 0 & 0 \\
          \frac{1}{2} & 1 & 1 & 0 \\
          \frac{1}{6} & \frac{1}{2} & 1 & 1
        \end{pmatrix},\quad
\end{equation}
as easily follows from $\ch(L \otimes \mathcal{E}) = \ch(L) \cup
\ch(\mathcal{E}) = e^H \cup \ch(\mathcal{E})$.

The twist $T_{\mathcal{O}_X}$ on the level of cohomology is given by
$\gamma \mapsto \gamma -\int_X\gamma \cup td(X) \cdot 1$ and hence its
matrix is given by
\begin{equation}
\label{eq:DMcon} 
  S = \begin{pmatrix}
          1 & -c & 0 & -d \\
          0 & 1 & 0 & 0 \\
          0 & 0 & 1 & 0 \\
          0 & 0 & 0 & 1
        \end{pmatrix}
\end{equation}
where
\[
  d:=H^3,\qquad c:=c_2\cdot H/12.
\]
The matrix $Q$ representing the bilinear form $\langle \cdot, \cdot
\rangle$ in this basis is given by
\[
  Q = \begin{pmatrix}
         0 & c  &   0  & d \\
         -c & 0 & -d &  0  \\
               0       &     d     &   0  &  0  \\
             -d      &      0      &   0  &  0
      \end{pmatrix}
\]
Now Kontsevich \cite{Ko2} observed the miracle that for the quintic
and its mirror the matrices $T$ and $S$ indeed correspond to
monodromy matrices of the Picard-Fuchs operator
\[
  \theta^4 - 5^5 z \bigl (\theta+\tfrac{1}{5}\bigr ) 
  \bigl (\theta+\tfrac{2}{5}\bigr ) \bigl (\theta+\tfrac{3}{5}\bigr )
  \bigl (\theta+\tfrac{4}{5} \bigr ).
\]
It has $0$, $1/5^5$ and $\infty$ as singular points. In an appropriate base,
the monodromy around $0$ is given by $T$ and around $1/5^5$ by $S$. 

We see that apparently the following happens: there is a point of
maximal unipotent monodromy, corresponding to $\otimes
\mathcal{O}(H)$ in $\Auteq(D^b(X))$ and there is a conifold point,
corresponding to the twist along $\mathcal{O}_X$.

Similar things occur in all the 14 hypergeometric cases. As there are
only three singular points in these cases, these two monodromies generate
the monodromy group. We refer to \cite{H} for a generalisation to
Calabi-Yaus in more general toric manifolds.

Calabi-Yau spaces with Picard number one seem to be rather scarse.  Apart
from the 14 hypergeometric cases there is there is a list (not
claiming completeness in any sense) by Borcea \cite{Bo} containing $11$
further cases. The examples are ramified covers and complete
intersections in Fano-varieties with Picard-number one. We know of a
few other cases. Basic invariants for such $X$ are the \emph{degree}
$d:=H^3$, the \emph{second Chern class} $c_2 \cdot H$ and the Euler
number $c_3=\chi_{top}$, of which the first two can be read off from
the matrix $S$.

It is sometimes more convenient to work with a different
representation based on the one used by C.~Doran and J.~Morgan (see
\cite{DM}).  That basis can be obtained from the one above using the
coordinate transformation given by the matrix
\[
  W = \begin{pmatrix}
        0 & 0 & 0 &  1 \\
        0 & 0 & 1 & -1 \\
       0 & \frac{1}{d} & -\frac{1}{2} &
          \frac{1}{3}-\frac{ c }{ d} \\
        \frac{1}{d} & 0 & -\frac{c }{ d} &
          \frac{c}{ d}
      \end{pmatrix}
\]
where $c$ and $d$ are as above.  This yields the following
representation:
\begin{gather*}
  T_{\text{DM}} = W^{-1} T W 
  = \begin{pmatrix}
      1 & 1 & 0 & 0 \\
      0 & 1 & d & 0 \\
      0 & 0 & 1 & 1 \\
      0 & 0 & 0 & 1
    \end{pmatrix},\quad
  S_{\text{DM}} = W^{-1} S W
  = \begin{pmatrix}
       1 & 0 & 0 & 0 \\
      -k & 1 & 0 & 0 \\
      -1 & 0 & 1 & 0 \\
      -1 & 0 & 0 & 1
    \end{pmatrix} \\
  Q_{\text{DM}} = W^t Q W
  = \begin{pmatrix}
       0 &  0 &  0 &  1 \\
       0 &  0 & -1 &  1 \\
       0 &  1 &  0 & -k \\
      -1 & -1 &  k &  0
    \end{pmatrix}
\end{gather*}
Here $d=H^3$ and $k=\frac{c_2 \cdot H}{12}+\frac{H^3}{6}$. This last
number has a simple interpretation as the dimension
$\textup{dim}(H^0(X,{\cal O}(H))$ of the linear system $|H|$.
 
\section{Computation of the monodromy}
Our starting point for the computation of the monodromy is the
following working hypothesis
\begin{hyp}
  Any differential equation of Calabi-Yau type which is strictly
  geometrical and for which the instanton numbers have an
  interpretation as the numbers of curves on a mirror manifold, the
  monodromy should satisfy the following conditions:
\begin{itemize}
\item[(H1)] There is a point of maximal unipotent monodromy,
  correspronding to $\otimes \mathcal{O}(H)$ in $\Auteq(D^b(X))$.
\item[(H2)] There is a conifold point, corresponding to the twist
  along $\mathcal{O}_X$.
\end{itemize}
\end{hyp}
By construction all the equations in the list from \cite{AESZ} have a
point of maximal unipotent monodromy at $z=0$. The non-obvious part is
to find a conifold point. We observed that in the cases where we know
the conifold point the \emph{spectrum}, i.e., the set of zeros of the
indicial equation at that point, was $\{0,1,1,2\}$. This is also
suggested by Hodge theory.  Therefore as a practical selection
criterium, we computed the indicial equations at the singular points
of all equations and found the equations with at least one singular
point with spectrum $\{0,1,1,2\}$.  As of the time of writing of this
article there were 178 such equations in our database. In many cases
there are several such points, but there are also some notable
exceptions, where no such singular point exists. An example is
equation \eqnno{32}, which is related to $\zeta(4)$ (see \cite{AZ}).
For the moment, we are unable to find integral or even just rational
lattices for these cases.

\begin{figure}
\begin{center}
\input{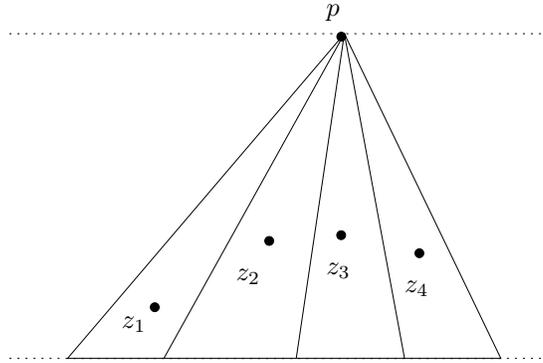}
\end{center}
\caption{Piecewise linear loops around the critical points $z_i$}
\label{fig:loops}
\end{figure}

For all the 178 equations that do have at least one singular point
with spectrum $\{0,1,1,2\}$ we computed high precision numerical
approximations for a set of generators of the monodromy group. These
computations were done in \texttt{Maple}. The first step was to
determine the critical points $z_1,\dots,z_\ell$ and to choose a
reference point $p$.  Next for each of the critical points $z_i$
except the point $z=\infty$ we choose a piecewise linear loop starting
and ending at the reference point $p$ and enclosing only one
critical point, namely $z_i$ (see \sref{fig:loops}).

Using the \texttt{Maple}-function \texttt{dsolve} we can numerically
integrate the differential equation along these paths. It turns out to
be a bit tricky to obtain the precision needed for the next steps. We
used the following options: method=gear, relerr=$10^{-15}$,
abserr=$10^{-15}$ and also increased \texttt{Digits} to 100. This
yielded the monodromy matrices with respect to an arbitrary basis and
produces fully filled $4 \times 4$-matrices with seemingly random
complex entries.

At this point there is a simple consistency check that we can do. If
there exists an integral lattice, the characteristic polynomial of
each of the monodromies should be a polynomial with integral
coefficients. As a further check, the roots of the indicial equations
at the corresponding singular points, should be logarithms of the
roots of the characteristic polynomial.  For the MUM-point and the
points with spectrum $\{0,1,1,2\}$ the characteristic polynomial
should be $(1-\lambda)^4$. This provides an indication of the
precision we have achieved.

The next step is to try and find a simultaneous base change that makes
all matrices integral, i.e., to find a monodromy invariant lattice
$\Lambda$. The crucial observation is the following. The monodromy $S$
around an $A_1$-singularity has the property that $\rk(S-\Id)=1$. The
one-dimensional image of $S-\Id$ is the span of the vanishing cycle.
Now choose one of the singular points with spectrum $\{0,1,1,2\}$ and
call the monodromy around the loop enclosing this singular point $S$.
As we are working with numerical approximations we cannot expect
$S-\Id$ to have rank $1$, but we can hope that the columns of the
matrix $S-\Id$ are nearly proportional. In that case we can pick an
arbitrary vector and apply $S-\Id$ to it. In this way we find a vector
$v_0$ that should be a good approximation to a lattice vector.

Further lattice vectors $v_1,\dots,v_k$ can be obtained by applying
words in the numerically computed monodromy matrices to $v_0$.  By
picking $n$ independent vectors among the ones found in this way, we
should find a basis for $\Lambda \otimes \Q$. When we transform the
monodromy matrices to this basis, the resulting matrices should have
rational entries. Of course this will not be exact, but we can try to
find rational matrices close to the matrices that we do find.  For
this we used continued fractions.  It may happen that we get very
large denominators or that the rational approximation is not very
accurate.  In that case we can try another set of $n$ independent
vectors among the $v_i$. If that is not successful, we can try another
point with spectrum $\{0,1,1,2\}$, if there is any. As a consistency
check, we can compute the characteristic polynomials of these rational
matrices and check that they have integral coefficients. As noticed
above, at the MUM-point and the conifold point the characteristic
polynomial should be $(1-\lambda)^4$, which we can also check. If any
of these checks fails, we have to try again with a different basis or
a different singular point with spectrum $\{0,1,1,2\}$. However, it
can and does happen that we try all potential conifold points and
several choices of a basis in each case, but do not find a rational
basis. We did find a rational basis in 143 of the 178 investigated
cases.

The rational basis found in this way is still rather arbitrary.
However, a major advantage is that at this point we expect to be
working with the \emph{exact} monodromy matrices. This allows us to do
linear algebra without worrying about the extra complications of
working with non exact numerical approximations. Provided that the
monodromy matrices around the MUM-point and the conifold point have
the right Jordan structure, we can find a new basis such that with
respect to this basis they have the standard form \pref{eq:DM0} and
\pref{eq:DMcon}. In a geometrical situation we expect the transformed
matrices to be integral. This happens in 64 cases. When we have the
monodromies around the MUM-point and the conifold point in the
standard form, we can read off the invariants $H^3$ and $c_2\cdot H$
and try to match the invariants with those of known Calabi-Yau spaces.

Despite our efforts to identify equivalent Calabi-Yau equations our
list probably still includes some Calabi-Yau equations that correspond
to the same geometrical situation. Transformations in the parameter
$z$ are a way of constructing seemingly different equations that
actually describe the same geometrical situation. In a geometrical
language this corresponds to pullback under a map $f: \P^1 \rightarrow
\P^1$. If the map $f$ is not injective, this may increase the number
of singular points. As long as the map $f$ is unramified around the
MUM-point and the conifold point it does not change the monodromies
and we ought to find the same $H^3$ and $c_2 \cdot H$. So as practical
way of trying to group together the equations that correspond to the
same geometry, we sort the 64 integral equations we found according to
$H^3$ and $c_2 \cdot H$. If we find several equations with the same
$H^3$ and $c_2 \cdot H$, it turns out that the genus zero instanton
numbers also coincide. This is a strong indication that these equations
are equivalent.

\section{Conifold-period and Euler characteristic}
If $\Omega$ is a family of holomorphic three forms on $Y_s$
(that is, a section of ${\cal L}:=\rho_*(\omega_{\cal Y}/S)$)
and $\Gamma$ is a horizontal family of cycles, then the
\emph{periods}
\[ 
  \int_{\Gamma} \Omega 
\]
are the solutions of the associated Picard-Fuchs equation.  In our
situation we identified two cycles, namely the torus $\mathbf{T}$ near
the MUM-point, and the vanishing sphere $\mathbf{S}$ near the conifold
point $z_c$. Correspondingly we have the fundamental period
$\int_{\mathbf{T}} \Omega$ which is the unique holomorphic solution
near the MUM-point. Equally important is the period $\int_{\mathbf{S}}
\Omega$ which we call the \emph{conifold-period} and which was called
$z_2(t)$ in the paper \cite{CdOGP}.  As the local monodromy around the
conifold point is supposed to be a symplectic reflection in $S$, this
special period can be determined directly from the differential
equation as follows. At such a conifold point there exists a basis of
solutions to the Calabi-Yau-equation around this point that consists
of three power series solutions and one solution of the form
\[
  y(z) = f(z) \log(z-z_c) + g(z).
\]
Going around $z_c$ once $y$ is replaced by $y(z) + 2\pi i f(z)$.  The
power series $f(z)$ represents a special solution around $z=z_c$ that
is determined up to a multiplicative scalar and which we call the
\emph{conifold-period}. This function can be continued analytically
around an arbitrary path which avoids the singularities of the
differential equation.  It is a remarkable fact that in all (but one,
namely nr.~\eqnno{224}) of the examples we know, the point $z_c$ is the
singular point that is \emph{closest to the origin}.  So there is a
preferred path from $z_c$ to $0$ by going along a straight line and we
consider the analytic continuation along this path. In \cite{CdOGP}
the expansion of the conifold period around $0$ is derived for the
quintic. It has the form
\begin{equation}
\label{eq:z2t}
  z_2(t) = \frac{H^3}{6} t^3 + \frac{c_2 \cdot H}{24} t +
  \frac{c_3}{(2\pi i)^3} \zeta(3) + O(q).
\end{equation}
The term $O(q)$ stands for any terms containing $q=e^{2\pi i t}$ and
$t = \frac{1}{2\pi i} \frac{y_1(z)}{y_0(z)} = \frac{1}{2\pi i} \log z
+ \frac{1}{2\pi i} \frac{f_1(z)}{f_0(z)}$ instead of \pref{eq:tcoord}.
Remarkable here is the `constant term' $\frac{c_3}{(2\pi i)^3}
\zeta(3)$.  This term is related to the four-loop correction to
the free energy $F_0$ introduced in \cite{CdOGP}.\footnote{To be
  precise, one has an expansion (see \cite{KKV}):
\[
  F_0=H^3 \frac{t^3}{3!} + (c_2 \cdot H)t + \frac{\chi}{2}\zeta(3) + 
    \sum_{d=1}^{\infty} n_d^0 Li_3(q^d)
\]
where $\mathrm{Li}_3(x):=\sum_{k=1}^{\infty} k^{-3} x^k$ is the classical
trilogarithm.}

One can conjecture this expansion to hold in all cases, which leads to
the following algorithm to determine $c_3$.  One can easily compute an
expansion of $y(z)$ to an arbitrary number of terms, e.g., using the
\texttt{Maple}-function \texttt{formal\_sol} from the
\texttt{DEtools}-package, as we did. This allows us to find $f(z)$ as
the coefficient of $\log(z-z_c)$. Around the MUM-point $z=0$ we can
compute expansions of the elements $y_i(z)$ of the Frobenius basis
(see \pref{app:inst}). We suppose that the domains of convergence of
the solutions around $z=0$ and those around $z=z_c$ overlap. That
enables us to pick some point $z_*$ where both expansions converge.
Computing numerically $f^{(k)}(z_*)$ ($k=0,\dots,3$) and
$y^{(k)}_i(z_*)$ ($k,i=0,\dots,3$), we can consider the equations
\[
  f^{(k)}(z_*) = \sum_{i=0}^3 c_i y_i^{(k)}(z_*).
\]
These equations can be solved for the $c_i$ and determine the analytic
continuation $z_2$ around $0$ of $f(z)$ as a linear combination of the
$y_i(z)$
\[
  z_2 = \sum_{i=0}^3 c_i y_i(z).
\]
From this we can readily read of the expansion of $z_2$ in $t$.  At
this point we can already check that the coefficient of $t^2$
vanishes.  As the conifold-period $f(z)$ was only determined up to a
constant, of course $z_2(t)$ ist determined up to a constant. If we
suppose that $H^3$ is known, then one can multiply the expansion for
$z_2(t)$ by a constant such that the coefficient of $t^3$ is
$\frac{H^3}{6}$.  We can then read off $c_2 \cdot H$ and $c_3$. In
praxis we find $H^3$ as discussed above from the monodromy generators.
This also yields $c_2 \cdot H$, so we have one more consistency check.
It is remarkable that in all cases we indeed find an integral value of
$c_3$!

\begin{table}
\vspace{-3cm}

\caption{\textbf{Calabi-Yau equations with integral monodromy}}
\begin{center}
\makebox[\textwidth]{%
\newcommand{\st}{\makebox[0mm][l]{*}}
\newcommand{\qm}{\makebox[0mm][l]{?}}
\begin{tabular}{|R|R|R|R||l|l|L|l|}
\hline
H^3 & c_2\cdot H & c_3 & |H| & Sings & Database &
\multicolumn{1}{l|}{Description} & Reference\\ \hline
1&10&  48\qm    & 1 &4* &\eqnno{225}&                                     &\\
1&22&-120    & 2 &3  &\eqnno{13}&X(6,6) \subset \P^5(1,1,2,2,3,3)     &\cite{KT}\\
1&34&-288    & 3 &3  &\eqnno{2} &X(10)   \subset \P^4(1,1,1,2,5)       &\cite{M1}\\
1&46& -484   & 4 &3  &\eqnno{9} &X(2,12) \subset \P^5(1,1,1,1,4,6) (?) &\cite{Al},\cite{DM}\\
2&20&  -44   & 2 &4* &\eqnno{271}&                                     &\\
2&32&-156    & 3 &3  &\eqnno{12}&X(3,4) \subset \P^5(1,1,1,1,1,2)      &\cite{KT}\\
2&44&-296    & 4 &3  &\eqnno{7}&X(8) \subset \P^5(1,1,1,1,4)           &\cite{M1}\\
3&42&-204    & 4 &3  &\eqnno{8}, \eqnno{125}&X(6) \subset \P^4(1,1,1,1,2) &\cite{M1}\\
4&40&-144    & 4 &3  &\eqnno{10}&X(4,4) \subset \P^5(1,1,1,1,2,2) &\cite{KT}\\
4&52&-256    & 5 &3  &\eqnno{14}, \eqnno{85}, \eqnno{86}&X(2,6)
\subset \P^5(1,1,1,1,1,3) &\cite{KT}\\
5&38& -100\st & 4 & 4* &\eqnno{302} & & \\
5&50&-200    & 5 &3  &\eqnno{1}, \eqnno{79}, \eqnno{87}, \eqnno{128}&X(5) \subset \P^4 &\cite{CdOGP}\\
5&62& -310 & 6 &4  &\eqnno{63}&                                      &        \\
6&36& -72   & 4 &4* &\eqnno{33}&                                      & \\
6&48&-156    & 5 &3  &\eqnno{11}, \eqnno{95}&X(4,6) \subset \P^5(1,1,1,2,2,3)      &\cite{KT}\\
7&46& -120\st & 5 &4* &\eqnno{109}&           &          \\
8&32& -8 & 4  & 4 &\eqnno{291}&           &          \\
8&56&-176    & 6 &3  &\eqnno{6}, \eqnno{75}, \eqnno{76}, \eqnno{96}&X(2,4) \subset \P^6  &\cite{LT}\\
9&30& 12\qm & 4 &4  &\eqnno{73}&                                     &         \\
9&54&-144    & 6 &3  &\eqnno{4}&X(3,3) \subset \P^5                    &\cite{LT}\\
10&40& -50   & 5 &5* &\eqnno{118}& & \\
10&40& -32 & 5 & 4* & \eqnno{292} & & \\
10&52& -116\st & 6 &4* &\eqnno{263}&                                     &\\ 
10&64&-200   & 7 &4  &\eqnno{51}&\text{\bfseries Conj: } X \stackrel{2:1}{\lra} B_5           &\cite[nr.\ 14]{Bo}\\
12&36& -32   & 5 &5* &\eqnno{117}&                                     &\\
12&48& -60   & 6 &5* &\eqnno{267}&                                     &\\
12&60&-144   & 7 &3  &\eqnno{5}, \eqnno{90}, \eqnno{91}, \eqnno{93}&X(2,2,3) \subset \P^6                  &\cite{LT}\\
13&58&-120   & 7 &4* &\eqnno{99}&\text{\bfseries Conj: } \text{$5\times5$-Pfaffian} \subset \P^6 &\cite{To}\\
14&56&-98    & 7 &5* &\eqnno{222}&\text{$7\times7$-Pfaffian} \subset
\P^6 &\cite{Ro98}\\
14&56& -100  & 7 & 5* & \eqnno{289} & & \\
15&54&-78    & 7 &?  &?  & \text{To}_{15} \subset \P^6          &\cite{To}\\
15&66&-150   & 8 &4  &\eqnno{24}&X(1,1,3) \subset \Grass(2,5)          &\cite{BCKvS}\\
16&52&-60    & 7 &?  &?  & \text{To}_{16} \subset \P^6          &\cite{To}\\
16&64&-128   & 8 &3  &\eqnno{3}, \eqnno{72}, \eqnno{224}&X(2,2,2,2) \subset \P^7        &\cite{LT}\\
17&50&-44    & 7 &?  &?  &\text{To}_{17}  \subset \P^6          &\cite{To}\\
18&60&-88    & 8 &4  &\eqnno{266}&                             &\\
20&68&-120   & 9 &4  &\eqnno{25}&X(1,2,2) \subset \Grass(2,5)  &\cite{BCKvS}\\
21&66&-102   & 9 &5* &\eqnno{254}&                             &\\
21&66&-100   & 9 &5* &\eqnno{270}&                             &\\
24&72&-116   & 10&4  &\eqnno{29}&\text{\bfseries Conj: } X(1,1,1,1,1,1,2) \subset X_{10} &\cite[nr.\ 6]{Bo}\\
25&70& -100\st & 10&5* &\eqnno{101}&                                       &\\
28&76&-116   & 11&4  &\eqnno{26}&X(1,1,1,1,2) \subset \Grass(2,6)        &\\
29&74& -100\st & 11&5* &\eqnno{256}&                                       &\\
32&80&-116   & 12&4  &\eqnno{42}&\text{\bfseries Conj: } X(1,1,2) \subset \LGrass(3,6)&\cite[nr.\ 8]{Bo}\\
33&78& -102\st & 12&5* &\eqnno{259}&                                       &\\
34&76& -88  & 12&4* &\eqnno{255}&                                       &\\
36&72& -72  & 12&5* &\eqnno{100}&                                       &\\
36&84&-120  & 13&4  &\eqnno{184}&\text{\bfseries Conj: } X(1,2) \subset X_5 &\cite[nr.\ 9]{Bo}\\
42&84&-98   & 14&6  &\eqnno{27}& X(1,1,1,1,1,1,1) \subset \Grass(2,7) &\cite{BCKvS}\\
42&84&-96   & 14&4  &\eqnno{28}&X(1,1,1,1,1,1) \subset \Grass(3,6)    &\cite{BCKvS} \\
44&92&-128  & 15&?  &?& X \stackrel{2:1}{\lra} {\text{$A_{22}$ or $A_{22}'$}}&\cite[nr.\ 10]{Bo}\\
47&86& -90\st & 15&6**&\eqnno{257}&                                     &\\
56&92& -92  & 17&?  &?&X(1,1,1,1) \subset F_1(Q_5) &\cite[nr.\ 24]{Bo}\\
57&90& -84  & 17&5* &\eqnno{247}& \text{Tj{\o}tta's example}             &\cite{Tj97}\\
\hline
\end{tabular}}
\end{center}
\label{tab:moni}
\end{table}

\section{Comments on the table of Calabi-Yau-equations}
In \sref{tab:moni} the heading Sings denotes the number of singular
point of the (first mentioned) differential equation. An additional
$*$ indicates, that an apparent singularity is present, around which
there is no monodromy.  The notation $X(\dots)$ denotes a complete
intersection of the indicated degrees in the indicated manifold.
Apart from the familiar $13$ hypergeometric cases and the cases from
complete intersection in Grassmanians that were studied in in
\cite{BCKvS}, one finds a few notable further cases.  First there is
the elusive $14$th hypergeometric case, observed in \cite{Al} and
\cite{DM}.  Any complete intersection $X(2,12)$ inside
$\P(1,1,1,1,4,6)$ has a singular point of type $A_1/(\Z/2)$, which
does not admit a crepant resolution. The case $X \stackrel{2:1}{\lra}
B_5$ is the Calabi-Yau double cover of the Fano-threefold $B_5$, which
is nothing but the three-dimensional section of $\Grass(2,5)$, which
is no.~14 in the list of Borcea. We found a fit with the equation
\eqnno{51} from \cite{AESZ}. A mirror for this Calabi-Yau is not
known, but we conjecture the Picard-Fuchs equation to be the indicated
one. We find similar fits for
\begin{description}
\item[\mathversion{bold}$X(1,1,1,1,1,1,2) \subset X_{10}$:] Here
  $X_{10} \subset \P^{15}$ is the celebrated 10-dimensional spinor
  variety of isotropic $4$-planes in the 8-dimensional quadric. 
\item[\mathversion{bold}$X(1,1,2) \subset \LGrass(3,6)$:]
  $\LGrass(3,6)=\mathrm{Sp}(3,\C)/P(\alpha_3) \subset \Grass(3,6)$ is
  the Lagrangian Grassmanian. 
\item[\mathversion{bold}$X(1,2) \subset X_5$:] Here
  $X_5=G_2/P(\alpha_{\text{long}}) \subset \Grass(5,7)$ is the space of
  5-di\-men\-sio\-nal subspaces isotropic for a 4-form on a 7-dimensional
  space. 
\end{description}
These are complete intersections inside homogeneous spaces. In
principle one can calculate the Picard-Fuchs equation for the
instanton numbers for these cases and verify our conjecture. The first
method consist in computing the quantum cohomology of these
homogeneous examples (for example by fixed point localisation) and then
use the quantum Lefschetz hyperplane principle. A second method
consists of finding a toric degeneration and then using polar duality. 
Such toric degenerations have been constructed for all spherical
varieties in \cite{AB}. Both methods were used in \cite{BCKvS} for the
case of complete intersections in Grassmannians. 

In his thesis \cite{To}, F. Tonoli considers Calabi-Yau varieties in
$\P^6$ of degree $12$ up to $17$. The first one is the complete
intersection $X(2,2,3)$, the second one the $5 \times 5$-Pfaffian, for
which we found a fit with the data from equation \eqnno{99}. The $7
\times 7$-Pfaffian was considered in \cite{Ro98}. The remaining three
case are new Calabi-Yau threefolds for which we have not yet found
corresponding Picard-Fuchs equations.

The column for the Euler characteristic $c_3$ was determined using the
expansion of the conifold-period around the MUM-point. It is a miracle
that we found integral values in all cases (except \eqnno{224}). This
checked with the known Euler number in those cases where a geometrical
interpretation was known. However, there are two notable cases where
we get a \emph{positive} value for $c_3$, which excludes an
interpretation as a Calabi-Yau space with Picard number one.
Furthermore, the conjectural integrality of elliptic intanton numbers
implies some congruence property on $c_3$. In most cases this was
satisfied, giving a strong indication that a Calabi-Yau threefold with
the indicated invariants should exist.  In some cases however, we
found non-integral in this way $n_d^1$. This is indicated with a $*$
after the value for $c_3$.

In the database column we indicate the number of the equation in the
electronic database of Calabi-Yau equations that can be found at the
web address
\[
 \text{\hlink{http://enriques.mathematik.uni-mainz.de/enckevort/db}}
\]
Up to 180 these numbers coincide with the ones used in \cite{AESZ}. 
For higher numbers one should check the source field in the electronic
database. If it contains Almkvist[$n$] the corresponding number in
\cite{AESZ} is $n$.

\section{Some Examples}
Let us now discuss a few typical examples from \sref{tab:moni} in more
detail. For full information on the other cases, we refer to the database
mentioned above.

\subsubsection*{Example 1}
The first equation we want to study is equation \eqnno{28} from
\cite{AESZ}, which is given by the following operator
\[
  L = \theta^4 - z(65 \,\theta^4+130 \,\theta^3 +105 \,\theta^2 
    + 40 \,\theta +6) + 4 z^2 (4\,\theta+3) (\theta+1)^2 (4\,\theta+5),
\]
where $\theta = z \frac{d}{dz}$.  This differential operator has four
singular points, namely $0$, $1/64$, $1$, and $\infty$. The Riemann
scheme is
\[
  P \left \{ \begin{array}{cccc}
      0& 1/64& 1& \infty\\ \hline 
      0&  0  & 0&  3/4\\
      0&  1  & 1&   1\\
      0&  1  & 1&   1\\
      0&  2  & 2&  5/4
    \end{array} \right \}
\]
Here the columns are the spectra, i.e., the set of solutions to the
indicial equation at the singular point indicated above the line.  The
points $1/64$ and $1$ have spectrum $\{0,1,1,2\}$, so they are
potential conifold points.  Using the algorithm discussed above we
computed the monodromies around the critical points and found an
integral lattice. With respect to this lattice the monodromy matrices
are as follows
\[
\begin{split}
  T&=T_0=\begin{pmatrix}
        1 & 1 &  0 & 0 \\
        0 & 1 & 42 & 0 \\
        0 & 0 &  1 & 1 \\
        0 & 0 &  0 & 1
      \end{pmatrix},\quad
  S=T_{\frac{1}{64}}=\begin{pmatrix}
        1  & 0 & 0 & 0 \\
       -14 & 1 & 0 & 0 \\
       -1  & 0 & 1 & 0 \\
       -1  & 0 & 0 & 1
      \end{pmatrix},\\
  T_1&=\begin{pmatrix}
          37 &  12 & -252 &  156 \\
        -126 & -41 &  882 & -546 \\
         -12 &	-4 &   85 &  -52 \\
         -18 &	-6 &  126 &  -77
      \end{pmatrix},\\
  T_{\infty}&=(T_1 T_{\frac{1}{64}} T_0)^{-1} =\begin{pmatrix}
         77 &  29 & -588 &  348\\
       -112 & -41 &  840 & -504\\
         -6 &  -2 &   43 &  -27\\
        -17 &  -6 &  126 &  -77
      \end{pmatrix}.
\end{split}
\]
From this, one can read off the invariants
\[
  H^3=42, \quad c_2\cdot H =84.
\]
In this case we know that the equation is the Picard-Fuchs equation
of the complete intersection $X(1,1,1,1,1,1)$ in $\Grass(3,6)$ and we
can easily check that these numbers coincide with the ones computed
from the geometry. Of course, the value $c_3$ computed from the
expansion of the conifold-period gives the right value $-96$.

One can easily check that $T_{\frac{1}{64}}$ and $T_1$ are of the
Picard-Lefschetz form $S_{\lambda,v}$, with the vector $v$ given by
\[
  v_{\frac{1}{64}}=\begin{pmatrix} 0 \\ 14 \\ 1 \\ 1
  \end{pmatrix},\quad
  v_1=\begin{pmatrix} 6 \\ -21 \\ -2 \\ -3 \end{pmatrix}.
\]
For $T_{\frac{1}{64}}$ we have $\lambda=1$, but for $T_1$ we have
$\lambda=2$. So the critical point $z=1$ is \emph{not} an ordinary
conifold point. This $\lambda=2$ is exactly what is needed to get
integral genus one instanton numbers with the recipe from
\pref{app:inst}. The first few elliptic instanton numbers that we find
in this way are
\[
  n^1_1=n^1_2=n^1_3=n^1_4=0,\quad n^1_5=84,\quad n^1_6=74382,\quad
  n^1_7=8161452.
\]
So it appears that there is a $\R P^3$ vanishing at the point $1$.
The derived category of coherent sheaves in a Grassmannian is
reasonably well understood (see \cite{Ka}) and so one can hope to
study in detail what happens in $D^b(X)$. This will be persued at
another place, \cite{EvS}.

 
\subsubsection*{Example 2}
Our second example has been discussed in \cite{Ro98,Tj99}. It is
interesting because there are two points with maximal unipotent
monodromy both of which have a geometrical interpretation. Because our
convention is to have the point of maximal unipotent monodromy that we
are considering at $z=0$ this example occurs twice in our list: once
as \eqnno{27} and once as \eqnno{222}.

In the former case the differential operator is given by
\[
\begin{split}
  L&=3^2\,\theta^4 -
  3z(173\,\theta^4+340\,\theta^3+272\,\theta^2+102\,\theta+15) \\
  &\quad -2z^2(1129\,\theta^4+5032\,\theta^3+7597\,\theta^2+4773\,\theta+1083) \\
  &\quad +2 z^3 (843\,\theta^4+2628\,\theta^3+2353\,\theta^2+675\,\theta+6) \\
  &\quad -z^4 (295\,\theta^4+608\,\theta^3+478\,\theta^2+174\,\theta+26) 
   +z^5 (\theta+1)^4
\end{split}
\]
The Riemann scheme of equation \eqnno{27} is
\[
  P \left \{ \begin{array}{cccccc}
      \zeta_1 & 0 & \zeta_2 & 3 & \zeta_3 & \infty\\ \hline 
         0    & 0 &    0    & 0 &    0    &   1\\
         1    & 0 &    1    & 1 &    1    &   1\\
         1    & 0 &    1    & 3 &    1    &   1\\
         2    & 0 &    2    & 4 &    2    &   1
    \end{array} \right \},
\]
where $\zeta_1<\zeta_2<\zeta_3$ are the (real) roots of
$z^3-289z^2-57z+1$.  The monodromies can be determined with our usual
recipe
\[
\begin{split}
  T_{\zeta_1}&=\begin{pmatrix}
    15 &   7 & -98 &  49 \\
   -28 & -13 & 196 & -98 \\
    -2 &  -1 &  15 &  -7 \\
    -4 &  -2 &  28 & -13
  \end{pmatrix},\quad
  T=T_0=\begin{pmatrix}
    1 & 1 &  0 & 0 \\
    0 & 1 & 42 & 0 \\
    0 & 0 &  1 & 1 \\
    0 & 0 &  0 & 1
  \end{pmatrix},\\
  S&=T_{\zeta_2}=\begin{pmatrix}
      1 & 0 & 0 & 0 \\
    -14 & 1 & 0 & 0 \\
     -1 & 0 & 1 & 0 \\
     -1 & 0 & 0 & 1
  \end{pmatrix},\quad
  T_3=\Id,\quad T_{\zeta_3}=\begin{pmatrix}
      1 & 0 &   0 &    0 \\
    -84 & 1 & 392 & -392 \\
     -9 & 0 &  43 &  -42 \\
     -9 & 0 &  42 &  -41
  \end{pmatrix},\\
  T_\infty&=(T_{\zeta_3} T_3 T_{\zeta_2} T_0 T_{\zeta_1})^{-1}=\begin{pmatrix}
      85 &   6 & -448 &   399 \\
    -266 & -13 & 1330 & -1232 \\
     -26 &  -1 &  127 &  -120 \\
     -42 &  -2 &  210 &  -195
  \end{pmatrix}.
\end{split}
\]
Here the monodromy operators $T_{\zeta_1}$, $T_{\zeta_2}$, and
$T_{\zeta_3}$ can be written in the Picard-Lefschetz form $S_{1,v}$
with the vector $v$ given by
\[
  v_{\zeta_1}=\begin{pmatrix} 7 \\ -14 \\ -2 \\ -2
  \end{pmatrix},\quad
  v_{\zeta_2}=\begin{pmatrix} 0 \\ 14 \\ 1 \\ 1
  \end{pmatrix},\quad
  v_{\zeta_3}=\begin{pmatrix} 0 \\ 28 \\ 3 \\ 3 \end{pmatrix}.
\]
The operator for equation \eqnno{222} can be obtained by replacing
$y(z)$ by $w^{-1}y(w^{-1})$ where $w=1/z$. In this way one finds the
operator
\[
\begin{split}
  L&=\theta^4 - z(295\,\theta^4+572\,\theta^3+424\,\theta^2+138\,\theta+17) \\
  &\quad +2z^2(843\,\theta^4+744\,\theta^3-473\,\theta^2-481\,\theta-101) \\
  &\quad -2z^3(1129\,\theta^4-516\,\theta^3-725\,\theta^2-159\,\theta+4) \\
  &\quad -3z^4(173\,\theta^4+352\,\theta^3+290\,\theta^2+114\,\theta+18)
    +3^2z^5(\theta+1)^4.
\end{split}
\]
The Riemann scheme of \eqnno{222} also follows from that of \eqnno{27}
\[
  P \left \{ \begin{array}{cccccc}
      1/\zeta_1 & 0 & 1/\zeta_3 & 1/3 & 1/\zeta_2 & \infty\\ \hline 
          0     & 0 &     0     &  0  &     0     &   1\\
          1     & 0 &     1     &  1  &     1     &   1\\
          1     & 0 &     1     &  3  &     1     &   1\\
          2     & 0 &     2     &  4  &     2     &   1
    \end{array} \right \}.
\]
For the monodromies the relation is not so obvious. Doing the standard
computation we find
\[
\begin{split}
  T_{\zeta_1^{-1}}&=\begin{pmatrix}
    29 &  14 & -98 &  49 \\
   -56 & -27 & 196 & -98 \\
    -8 &  -4 &  29 & -14 \\
   -16 &  -8 &  56 & -27
  \end{pmatrix},\quad
  T=T_0=\begin{pmatrix}
    1 & 1 &  0 & 0 \\
    0 & 1 & 14 & 0 \\
    0 & 0 &  1 & 1 \\
    0 & 0 &  0 & 1
  \end{pmatrix},\\
  S&=T_{\zeta_3^{-1}}=\begin{pmatrix}
     1 & 0 & 0 & 0 \\
    -7 & 1 & 0 & 0 \\
    -1 & 0 & 1 & 0 \\
    -1 & 0 & 0 & 1
  \end{pmatrix},\quad
  T_{1/3}=\Id,\quad T_{\zeta_2^{-1}}=\begin{pmatrix}
       1 & 0 &   0 &    0 \\
    -105 & 1 & 294 & -294 \\
     -25 & 0 &  71 &  -70 \\
     -25 & 0 &  70 &  -69
  \end{pmatrix},\\
  T_\infty&=(T_{\zeta_2^{-1}} T_{1/3} T_{\zeta_3^{-1}} T_0
  T_{\zeta_1^{-1}})^{-1}=\begin{pmatrix}
     155 &  13 & -476 &   427 \\
    -420 & -27 & 1260 & -1162 \\
     -76 &  -4 &  225 &  -211 \\
    -126 &  -8 &  378 &  -349
  \end{pmatrix}.
\end{split}
\]
Again the monodromies $T_{\zeta_1^{-1}}$, $T_{\zeta_3^{-1}}$, and
$T_{\zeta_2^{-1}}$ can be written in the Picard-Lefschetz form
$S_{1,v}$ with the vector $v$ given by
\[
  v_{\zeta_1^{-1}}=\begin{pmatrix} 7 \\ -14 \\ -2 \\ -4
  \end{pmatrix},\quad
  v_{\zeta_3^{-1}}=\begin{pmatrix} 0 \\ 7 \\ 1 \\ 1
  \end{pmatrix},\quad
  v_{\zeta_2^{-1}}=\begin{pmatrix} 0 \\ 21 \\ 5 \\ 5 \end{pmatrix}.
\]
Despite the fact that we are really dealing with the same equation in
a different formulation, the monodromies look rather different. Of
course the monodromy groups generated by these matrices are
isomorphic, but it is not so easy to see.

\subsubsection*{Example 3}
The next example is equation \eqnno{29} from \cite{AESZ}. The operator
is
\[
  L = \theta^4 - 2z(2\,\theta+1)^2 (17\,\theta^2+17\,\theta+5) + 
      2^2 z^2 (2\,\theta+1)(\theta+1)^2(2\,\theta+3).
\]
In this case the Riemann scheme is
\[
  P \left \{ \begin{array}{cccc}
      0& \zeta_1 & \zeta_2 & \infty\\ \hline 
      0&    0    &    0    &  1/2\\
      0&    1    &    1    &   1\\
      0&    1    &    1    &   1\\
      0&    2    &    2    &  3/2
    \end{array} \right \},
\]
where $\zeta_1<\zeta_2$ are the (real) roots of $1-136z+16z^2$.  The
monodromy matrices are
\[
\begin{split}
T&=T_0=\begin{pmatrix}
  1 & 1 &  0 & 0 \\
  0 & 1 & 24 & 0 \\
  0 & 0 &  1 & 1 \\
  0 & 0 &  0 & 1
\end{pmatrix},\quad
S=T_{\zeta_1}=\begin{pmatrix}
    1 & 0 & 0 & 0 \\
  -10 & 1 & 0 & 0 \\
   -1 & 0 & 1 & 0 \\
   -1 & 0 & 0 & 1
\end{pmatrix},\\
T_{\zeta_2}&=\begin{pmatrix}
    51 &  20 & -240 &  140 \\
  -130 & -51 &  624 & -364 \\
   -15 &  -6 &   73 &  -42 \\
   -25 & -10 &  120 &  -69
\end{pmatrix},\\
T_\infty&=(T_{\zeta_2} T_{\zeta_1} T_0)^{-1} = \begin{pmatrix}
    71 &  31 & -360 &  200 \\
  -120 & -51 &  600 & -340 \\
   -10 &  -4 &   49 &  -29 \\
   -24 & -10 &  120 &  -69
\end{pmatrix}.
\end{split}
\]
One can check that the $T_{\zeta_i}$ can be written as $S_{1,v}$ with
$v$ given by
\[
  v_{\zeta_1}=\begin{pmatrix} 0 \\ 10 \\ 1 \\ 1
  \end{pmatrix},\quad
  v_{\zeta_2}=\begin{pmatrix} 10 \\ -26 \\ -3 \\ -5 \end{pmatrix}.
\]
From the expressions for $T$ and $S$ we find $H^3=24$, $c_2 \cdot
H=72$. We also have enough information to compute the elliptic
instanton numbers as a function of $c_3$. By equating $n^1_1=0$ we
find $c_3=-116$ and all the $n^1_i$ we computed are integral.  Thes
same value for $c_3$ is obtained from the expansion of the
conifold-period. It turns out that we are \emph{lucky} and that there
is exactly one 1-parameter Calabi-Yau known with these invariants,
namely $X(1,1,1,1,1,1,2) \subset X_{10}$ (see \cite{Bo}). So we
conjecture that equation \eqnno{29} is the Picard-Fuchs equation
corresponding to this Calabi-Yau. In the same way we conjecturally
identified the Picard-Fuchs equations of four more 1-parameter
Calabi-Yau spaces. We labelled these equations in \sref{tab:moni} by
writing \textbf{Conj:} in front of the conjectured Calabi-Yau.

\subsubsection*{Example 4}
As our final example we will use equation \eqnno{270} (218 in the
numbering from \cite{AESZ}) which is given by the diffential operator
\[
\begin{split}
  L &= 7^2\,\theta^4 -
  42z(192\,\theta^4+396\,\theta^3+303\,\theta^2+105\,\theta+14) \\
  &\quad +2^2 \cdot 3 z^2 (1188\,\theta^4+11736\,\theta^3 + 
     20431\,\theta^2+12152\,\theta+2436) \\
  &\quad +2^2 \cdot 3^3 z^3 (532\,\theta^4+504\,\theta^3 - 3455\,\theta^2 -
     3829\,\theta-1036) \\
  &\quad -6^4z^4(2\,\theta+1)(36\,\theta^3+306\,\theta^2+421\,\theta+156) \\
  &\quad -2^6\cdot 3^4 z^5 (2\,\theta+1)(3\,\theta+2)(3\,\theta+4)(2\,\theta+3).
\end{split}
\]
The Riemann scheme is
\[
  P \left \{ \begin{array}{cccccc}
      -7/12 & 0 & \zeta_1 & \zeta_2 & \zeta_3 & \infty\\ \hline 
        0   & 0 &    0    &    0    &    0    &   1/2\\
        1   & 0 &    1    &    1    &    1    &   2/3\\
        3   & 0 &    1    &    1    &    1    &   4/3\\
        4   & 0 &    2    &    2    &    2    &   3/2
    \end{array} \right \},
\]
where $\zeta_1$ is the real root of $1296z^3-864z^2+168z-1$ and
$\zeta_2$, $\zeta_3$ are its complex roots with $\im \zeta_2<0$ and
$\zeta_3=\bar{\zeta}_2$. The monodromies can be computed and turn out
to be integral
\[
\begin{split}
T_{-\frac{7}{12}}&=\Id,\quad T=T_0=\begin{pmatrix}
  1 & 1 &  0 & 0 \\
  0 & 1 & 21 & 0 \\
  0 & 0 &  1 & 1 \\
  0 & 0 &  0 & 1
\end{pmatrix},\quad T_{\zeta_1}=\begin{pmatrix}
   1 & 0 & 0 & 0 \\
  -9 & 1 & 0 & 0 \\
  -1 & 0 & 1 & 0 \\
  -1 & 0 & 0 & 1
\end{pmatrix},\\
T_{\zeta_2}&=\begin{pmatrix}
   16 &   5 & -60 &   40 \\
  -45 & -14 & 180 & -120 \\
   -6 &  -2 &  25 &  -16 \\
   -9 &  -3 &  36 &  -23
\end{pmatrix},\quad T_{\zeta_3}=\begin{pmatrix}
   11 &  5 & -45 &  25 \\
  -18 & -8 &  81 & -45 \\
   -2 & -1 &  10 &  -5 \\
   -4 & -2 &  18 &  -9
\end{pmatrix},\\
T_\infty&=(T_{\zeta_3} T_{\zeta_2} T_{\zeta_1} T_0 T_{-\frac{7}{12}})^{-1}=
\begin{pmatrix}
    8 &  5 & -45 &  21 \\
  -12 & -5 &  60 & -36 \\
   -1 &  0 &   4 &  -4 \\
   -3 & -1 &  15 & -10
\end{pmatrix}.
\end{split}
\]
So we find $H^3=21$ and $c_2 \cdot H=66$. The $T_{\zeta_i}$ can be
written in Picard-Lefschetz form with $\lambda=1$ and $v$ given by
\[
  v_{\zeta_1}=\begin{pmatrix} 0 \\ 9 \\ 1 \\ 1
  \end{pmatrix},\quad
  v_{\zeta_2}=\begin{pmatrix} 5 \\ -15 \\ -2 \\ -3
  \end{pmatrix},\quad
  v_{\zeta_3}=\begin{pmatrix} 5 \\ -9 \\ -1 \\ -2 \end{pmatrix}.
\]
We can also compute the elliptic instanton numbers. Setting $n^1_1=0$
we find $c_3=-100$ and with this value of $c_3$ all computed $n^1_d$
turn out to be integers. The same value of $c_3$ was obtained from the
expansion of the conifold-period. So we have a Calabi-Yau equation
that as far as we can check looks like the Picard-Fuchs equation of a
Calabi-Yau manifold. However, we do not know any 1-parameter
Calabi-Yau with the geometric invariants that we computed. In
\sref{tab:moni} there are some more equations which look geometrical
in every respect, but for which we have not found any Calabi-Yau yet.

\section{Open problems}
The work described in this paper is no more than a start and there are
many open problems left. We have found quite a few Calabi-Yau
equations that look in every respect like the Picard-Fuchs equation of
a Calabi-Yau manifold, but for which we do not know if a Calabi-Yau
manifold exists. We know the degree, the second Chern class, the Euler
characteristic and the instanton numbers.

To determine an integral lattice we need to single out two singular
points, where we bring the monodromies into the the standard forms
$T_{\text{DM}}$ and $S_{\text{DM}}$. When there are no singular points
with spectrum $\{0,1,1,2\}$ we do not have a good candidate for
$S_{\text{DM}}$ and cannot even start our procedure for determining an
integral lattice. It would be interesting to see what can be done in
such cases.  We also did not look for other integral lattices as in
\cite{DM}.

The conjectural appearance of the constant term $c_3\zeta(3)/(2\pi
i)^3$ in the expansion of the conifold period (and the free energy) is
very intriguing. Is this a mathematical theorem?

The key obstacle to computing the elliptic instanton numbers is
finding the holomorphic function $f(z)$ in \pref{eq:F1}. Our ansatz in
combination with our recipe for determining the exponents works a many
cases, but it is no more than an educated guess. A better understanding of 
the genus one computation in terms of the BCOV-torsion as in \cite{FL} 
will probably be helpful.

Many of the equations from the list in \cite{AESZ} come from Hadamard
products. The singular points of a Hadamard product are given by
products of the singular points of the factors. Maybe it is also
possible to determine the monodromies of the Hadamard product in terms
of the monodromies of the factors. 


\appendix
\section{Orbifolds of $A_1$}
In many examples one encounters monodromy tranformations that
are not described by the usual Picard-Lefschetz formula, but
rather are \emph{powers} of such operations. We offer 
a possible explanation of this phenomonen, which is only visible
on the integral level. 

Consider a lattice $\Lambda$ with bilinear form $\langle \cdot, \cdot
\rangle$. For $\beta \in \Lambda$ and $\lambda \in \Z$ consider the
the transformation
\begin{equation}
\label{eq:Slb}
  S_{\lambda,\beta}: \Lambda \lra \Lambda,\quad 
  S_{\lambda,\beta}(\alpha) = \alpha - \lambda \langle \beta, 
    \alpha \rangle \beta. 
\end{equation}
The transformation $S_{\lambda,\beta}$ preserves $\langle \cdot, \cdot
\rangle$ in the symmetric case only when $\lambda=2/Q(\beta,\beta)$
(or $\lambda=0$). In that case $S_{\lambda,\beta}$ has order two and is
a reflection.  When $\langle \cdot, \cdot \rangle$ is antisymmetric,
there is no restriction on $\lambda$ and $S_{\lambda,\beta} \circ
S_{\lambda',\beta} = S_{\lambda+\lambda',\beta}$. So in that case
$S_{\lambda,\beta}$ does not have finite order. 

Such transformations occur as monodromy transformations 
where not a sphere, but rather a quotient $S^3/G$ by a finite
group $G$ is vanishing, as we will explain now. 
Consider the function defining the three-dimensional $A_1$-singularity:
\[
  f: \C^4 \lra \C,\quad f(x,y,z,t)=x^2+y^2+z^2+t^2
\]
The fibre $F_s$ of $f$ over $s \in \C \setminus {0}$ is called the
Milnor fibre and can be identified with the cotangent bundel to the
sphere $\{(x,y,z,t) \in \R^4 \mid x^2+y^2+z^2+t^2=s\}$, which is
vanishing when $s \rightarrow 0$. We choose an orientation and let
$\delta$ be the homology class of this sphere. There is also a
covanishing cycle $\epsilon$ in the dual group
$\Hgrp^{\text{cl}}_3(F_s,\Z)$ (homology with closed support).  One has
\[
  \Hgrp_3(F_s,\Z) = \Z \delta,\quad 
  \Hgrp_3^{\text{cl}}(F_s,\Z) = \Z\epsilon,\quad 
  \langle \delta, \epsilon \rangle = 1
\]
Let $G \subset \mathrm{SU}(2) = S^3$ be a finite subgroup. $G$ then
acts linearly on $\R^4$ and by complexification on $\C^4$, leaving
invariant the function $f$ defining the $A_1$-singularity.  Consider
the quotient map $\pi: \C^4 \lra X := \C^4/G$. The space $X$ will be
singular, but $f$ descends to a function $g: X \lra \C$, such that
$f=\pi \circ g$. So the fibre $G_s := g^{-1}(s)$ is the quotient of
$F_s$ by $G$. In the fibre $G_s$ there is a cycle $S^3/G$ vanishing
when $s \rightarrow 0$, with homology class $d \in \Hgrp_3(G_s,\Z)$. 
As above there also exists a covanishing cycle $e \in
\Hgrp_3^{\text{cl}}(G_s,\Z)$ such that
\[
  \Hgrp_3(G_s,\Z) = \Z d,\quad \Hgrp_3^{\text{cl}}(G_s,\Z) = \Z
  e,\quad \langle d, e \rangle=1. 
\]
The map $\pi$ induces maps $\pi^*$ and $\pi_*$ between the homology
groups of $F_s$ and $G_s$ and one easily sees that
\[
  \pi_*(\delta) = |G|d,\quad \pi_*(\epsilon) = e,\quad 
  \pi^*(d) = \delta,\quad \pi^*(e)=|G|\epsilon
\] 
The Picard-Lefschetz formula tells us that under the monodromy of $f$
the cycle $\delta$ remains fixed, whereas the cycle $\epsilon$ gets
mapped to $\epsilon -\delta$. From the fact that the monodromy
commutes with the group action we obtain, by taking $\pi_*$, that $d$
remains fixed, whereas $e$ gets mapped to $e-|G|d$.  From this one
deduces in the usual way that the occurence of a singularity of type
$A_1/G$ will lead the monodromy transformation described by the
modified Picard-Lefschetz formula (see \cite{AGV})
\[
  \gamma \mapsto \gamma -|G|\langle d, \gamma \rangle d
\]  
The cycle $d$ should should give rise to a spherical object in the
derived category of the mirror, but only the $|G|$th power of the
Seidel-Thomas twist would arise from a monodromy transformation. 

\section{Computation of instanton numbers}
\label{app:inst}
According to \cite{GV1,GV2} (see also \cite{KKV}) we have the
following expansion for the partition function $F$ of the topological
string
\[
  F =\sum_{g=0}^\infty \lambda^{2g-2} F_g = \sum_{g=0}^\infty \sum_d
  \sum_{m=0}^\infty n^g_d \frac{1}{m} \biggl ( 2 \sin
  \frac{m\lambda}{2} \biggr )^{2g-2} q^{dm}. 
\]
The partion function $F$ can be defined physically or mathematically
using Gromov-Witten invariants. The above formula can then be
considered to define the Gopakumar-Vafa invariants $n^g_d$. In
contrast to e.g., the Gromov-Witten invariants, the Gopakumar-Vafa
invariants are conjectured to be always integral. 

We will restrict to the genus zero and genus one invariants. 
Furthermore, we will only consider the 1-parameter case. In that case
we have the following formulas (see \cite{KKV})
\[
  \partial_t^3 F_0 = n^0_0 + \sum_{\ell=1}^\infty \frac{n^0_\ell \ell^3
  q^\ell}{1-q^\ell}
\]
with $n^0_0=H^3$ and
\[
  \partial_t F_1 = \frac{c_2 \cdot H}{24} + \sum_{d=1}^\infty
  \sum_{k=1}^\infty \bigl ( \tfrac{1}{12} n_d^0 + n_d^1 \bigr ) d
  q^{kd}. 
\]
Here we define the coordinate $t$ by $q=e^{-t}$. So if we can compute the
left hand sides of these equations the invariants $n^0_d$ and $n^1_d$
can easily be determined. 

To do so, we first introduce a special basis of solutions for the
equation \pref{eq:L} around a point of maximal unipotent monodromy,
i.e., a singular point where $\lambda=0$ is the only solution to the
indicial equation. In physical terms such a point (also called
MUM-point) corresponds to a large radius limit point. 

Suppose $z=0$ is a MUM-point. Then we can use the Frobenius method. 
The idea is to consider a solution with values in the ring
$\C[\rho]/(\rho^n)$. We make the following ansatz for such a solution
\[
  \tilde{y}(z)=\sum_{n=0}^\infty A(n,\rho) z^{n+\rho}
  = y_0(z) + y_1(z) \rho + \dots + y_{n-1}(z) \rho^{n-1},
\]
where we define
\[
  z^\rho = e^{\log z \cdot \rho} = 1 + \log z \cdot \rho +
  \frac{\log^2 z}{\rho} \cdot \rho^2 + \dots + \frac{\log^{n-1}
    z}{(n-1)!} \cdot \rho^{n-1}. 
\]
Using $\theta z^{n+\rho} = (n+\rho) z^{n+\rho}$, where $\theta=z
\frac{d}{dz}$, we can translate the equation $L \tilde{y}=0$ into a
recursion relation for the $A(n,\rho)$. As initial condition for the
recursion we use $A(0,\rho)=1$. The $y_i$ we find in this way are
called the \emph{Frobenius basis}. 

Define power series $f_i$ by the following expression
\[
  \sum_{n=0}^\infty A(n,\rho) z^n = f_0(z) + f_1(z) \rho + \dots +
  f_{n-1}(z) \rho^{n-1}. 
\]
Because $z^\rho \sum_{i=0}^{n-1} f_i(z) \rho^i = \sum_{i=0}^{n-1}
y_i(z) \rho^i$, we find
\[
  y_i(z)=\sum_{j=0}^i \frac{\log^i z}{i!} f_{j-i}(z). 
\]
Using the Frobenius base we can define a new coordinate
\begin{equation}
\label{eq:tcoord}
  t=y_1(z)/y_0(z)=\log z + \frac{f_1(z)}{f_0(z)}. 
\end{equation}

There are basically two ways to compute $\partial_t^3 F_0$. The
starting point of the first one is the \emph{Yukawa coupling} in the
$z$ coordinate
\[
  K_{zzz}=\exp\bigl ({\textstyle -\frac{1}{2} \int a_3(z) \dop z} \bigr ),
\]
where $a_3(z)$ is one of the coefficients from \pref{eq:L}.  The claim
is that $\partial_t^3 F$ is the following transformation of this
function to the $t$-coordinate defined in \pref{eq:tcoord}
\begin{equation}
\label{eq:Kttt}
  \partial_t^3 F(t) = \frac{K_{zzz}(z(t))}{y_0^2(z(t))
    \bigl ( \frac{dt}{dz} \bigr )^3}
\end{equation}

Now recall that a Calabi-Yau equation has to satisfy a list of
conditions (see \cite{AZ,AESZ}). One of these can be written as
\begin{equation}
\label{eq:AZ}
  a_1=\frac{1}{2}a_2a_3-\frac{1}{8}a_3^3 + a_2' - \frac{3}{4}a_3a_3'-
  \frac{1}{2}a_3''
\end{equation}
According to Proposition~1 from \cite{AZ} this condition is equivalent
to the two conditions
\begin{align}
\label{eq:AZ1}
  \frac{d^2}{dt^2} \frac{y_2}{y_0} &= \frac{\exp\bigl (-\frac{1}{2}\int
    a_3(z) \dop z \bigr )}{y_0^2 \bigl ( \frac{dt}{dz}\bigr )^3}, \\
\label{eq:AZ2}
  \frac{d^2}{dt^2} \frac{y_3}{y_0} &= t \frac{d^2}{dt^2} \frac{y_2}{y_0}. 
\end{align}
The second condition is equivalent to the existence of a function $G$
and a constant $c$ such that
\[
  \Pi(t):=\frac{1}{y_0} \begin{pmatrix} y_0 \\ y_1 \\ y_2 \\ y_3
  \end{pmatrix} =
  \begin{pmatrix}
    1 \\
    t \\
    \partial_t G - c\\
    t\partial_t G - 2 G
  \end{pmatrix}. 
\]
The vector $\Pi(t)$ is called the \emph{normalized period vector}. 
Using the second condition \pref{eq:Kttt} translates to
\[
  \partial_t^3 F_0 = \partial_t^3 G(t) = \partial_t^2 \frac{y_2}{y_0}. 
\]
This yields a different way of computing $\partial_t^3 F_0$ and
therefore also the instanton numbers $n^0_d$. 

To compute $\partial_t F_1$ we use the recipe from \cite{BCOV1,BCOV2}
based on an analysis of the so called holomorphic anomaly. We will use
the following formula from \cite{BCOV1} (using our notation and
adapted slightly for the case we are studying):
\begin{equation}
\label{eq:F1}
  \partial_t F_1 = \partial_t \log \left (
    \frac{z^{1+\frac{c_2 \cdot H}{12}}f(z)}{y_0^{4-\frac{c_3}{12}}
    \frac{\partial t}{\partial z}} \right ). 
\end{equation}
In this formula one needs the geometrical data $c_2 \cdot H$ and $c_3$
which can usually be determined from the monodromy calculation and/or
conifold period. However, the main problem with this formula is the
function $f$ which is a holomorphic function of $z$ that still has to
be determined. We will use an ansatz for $f$ to reduce this problem to
the determination of a finite number of parameters. To describe this
ansatz note that because of the special form of a Calabi-Yau equation
we can write
\[
  a_4(z)= z^4 \Delta(z)= z^4 \prod_i (\Delta_i(z))^{k_i},
\]
for some polynomial $\Delta(z)$, which we call the
\emph{discriminant}.  The $\Delta_i(z)$ are the irreducible factors
(over $\R$) of $\Delta(z)$. Our ansatz is then the following
\[
  f(z) = \prod_i (\Delta_i(z))^{s_i},
\]
where the exponents $s_i \in \Q$ still have to be determined. The
(apparent) singular points of the operator are the zeros of the
discriminant $\Delta(z)$ (and $0$ and $\infty$). So each of the
factors $\Delta_i(z)$ corresponds via its zeros to one or more
(apparent) singular points. To determine the exponents we look at the
monodromies around the corresponding singular points. When the
singular point is a conifold, i.e., the monodromy is of the form
$S_{1,v}$, then the exponent is generally assumed to be
$-\frac{1}{6}$. We generalize this to $-\frac{\lambda}{6}$ for
monodromies of the form $S_{\lambda,v}$ for arbitrary $\lambda$. When
the monodromy is the identity, we put the exponent to zero. These
rules already allow us to deal with many equations. However,
monodromies of other types for which we do not know of a sensible
guess do occur.

\bibliographystyle{plain}
\bibliography{jointproj}
\end{document}